\documentclass[12pt]{article}
\usepackage{amsfonts}
\usepackage{a4wide}
\usepackage{amsmath,amscd,amsthm,a4,amssymb}
\usepackage{a4wide}

\begin{document}

\begin{center}
{\LARGE \textbf{A note on the approximation by a sum of two algebras}}

\

\textsc{\large Aida Kh. Asgarova}$^{1}$ \textsc{\large and} \textsc{\large %
Vugar E. Ismailov}\footnote{%
Corresponding author}$^{2}$ \vspace{3mm}

$^{1,2}${Institute of Mathematics and Mechanics, Baku, Azerbaijan} \vspace{%
1mm}

$^{1}${Department of Mathematics, Khazar University, Baku, Azerbaijan}
\vspace{1mm}

$^{2}${Center for Mathematics and its Applications, Khazar University, Baku,
Azerbaijan} \vspace{1mm}

e-mails: $^{1}$aidaasgarova@gmail.com, $^{2}$vugaris@mail.ru
\end{center}

\smallskip

\textbf{Abstract.} We consider the problem of approximation of a continuous
function $f$ defined on a compact metric space $X$ by elements from a sum of
two algebras. We prove a de la Vall\'{e}e Poussin type theorem, which
estimates the approximation error $E(f)$ from below. We also obtain a
duality formula for the precise computation of $E(f)$.

\smallskip

\textit{2010 MSC:} 41A30, 41A50, 46B50, 46E15

\smallskip

\textit{Keywords:} de la Vall\'{e}e Poussin theorem; the approximation
error; bolt; duality formula

\bigskip

\bigskip

\begin{center}
{\large \textbf{1. Introduction}}
\end{center}

The classical de la Vall\'{e}e Poussin \cite{Vallee} theorem says that if
the difference $f-P_{0}$ between a function $f\in C[a,b]$ and a polynomial $%
P_{0}$ of degree at most $n$ takes alternating signs at points $x_{i}\in
\lbrack a,b]$ of an ordered sequence of $n+2$ points, then

\begin{equation*}
\rho (f,\mathcal{P}_{n})\geq \min_{1\leq i\leq n+2}\left\vert
f(x_{i})-P_{0}(x_{i})\right\vert ,
\end{equation*}%
where $\mathcal{P}_{n}$ is the set of polynomials of degree at most $n$ and $%
\rho (f,\mathcal{P}_{n})$ is the distance from $f$ to $\mathcal{P}_{n}$, i.e.

\begin{equation*}
\rho (f,\mathcal{P}_{n})\overset{def}{=}\inf_{P\in \mathcal{P}%
_{n}}\max_{x\in \lbrack a,b]}\left\vert f(x)-P(x)\right\vert .
\end{equation*}%
The distance $\rho (f,\mathcal{P}_{n})$ is also called the approximation
error of $f$ by the class of polynomials $\mathcal{P}_{n}$.

In this paper we obtain an analogue of this result for approximation by a
sum of two algebras in the space of continuous functions. To make the
problem precise, assume $X$ is a compact metric space, $C(X)$ is the space
of real-valued continuous functions on $X$, $A_{1}$ and $A_{2}$ are closed
subalgebras of $C(X)$ that contain constants. For a given function $f\in
C(X) $ consider the approximation of $f$ by elements of $A_{1}+A_{2}$. We
are interested in estimating and/or precisely computing the approximation
error

\begin{equation*}
E(f)=\rho (f,A_{1}+A_{2})=\inf_{u\in A_{1}+A_{2}}\left\Vert f-u\right\Vert .%
\eqno(2.1)
\end{equation*}%
Here $\left\Vert \cdot \right\Vert $ denotes the standard max-norm in $C(X)$%
. If there exists a function $u_{0}\in A_{1}+A_{2}$ such that the $\inf $ in
(2.1) is attained $u_{0}$ is called a best approximation (or an extremal
element).

It should be remarked that numerous papers have addressed approximation
problems related to the sums of algebras (see Khavinson \cite{12} for
related results and references). The historical roots of this subject trace
back to M.H. Stone's papers \cite{St1,St2} published in 1937 and 1948, where
he specifically explored the case of approximation when only one algebra is
involved. The Stone-Weierstrass theorem, a noteworthy result in this
context, posits that a subalgebra $A\subset C(X)$, containing a nonzero
constant function, is dense in the entire space $C(X)$ if and only if $A$
separates points in $X$. Here, separation implies that for any distinct
points $x$ and $y$ in $X$, there exists a function $g\in A$ with $g(x)\neq
g(y)$. The density of the sum of two subalgebras $A_{1}$ and $A_{2}$ in $%
C(X) $ for a compact Hausdorff space $X$ has been extensively examined by
Marshall and O'Farrell \cite{Mar1,Mar2}. In \cite{Mar1} they obtained a
necessary and sufficient geometrical condition for the density of $%
A_{1}+A_{2}$ in $C(X)$. In \cite{Mar2} a description of Borel measures on $X$
orthogonal to the sum $A_{1}+A_{2}$ was provided. This description led to
another condition that is equivalent to the density of $A_{1}+A_{2}$ in $%
C(X) $. The open questions associated with the sum of more than two algebras
are also outlined in \cite{Mar2}.

Note that for many useful algebras $A_{1},A_{2}$ and compact spaces $X$, the
density property (or approximation with arbitrary precision) does not hold;
hence the problem of error analysis arises. Within which approximation error
we can approximate a given function $f\in C(X)$ by elements from $%
A_{1}+A_{2} $? The present paper addresses this natural question. Our work
is a continuation of the papers \cite{AI,AHI}, where the Chebyshev type
criterion for a best approximation from $A_{1}+A_{2}$ to $f$ was obtained.
In this paper, using the ideas from \cite{AI,AHI}, we prove a de la Vall\'{e}%
e Poussin type theorem, which estimates the approximation error $E(f)$ from
below. We also obtain a duality formula for the precise computation of $E(f)$%
.

The algebras considered in this paper are abstract. In particular cases,
they may turn into certain function algebras $A_{i}$, which are useful in
various applications. For example, these $A_{i}$ can be algebras of
univariate functions, ridge functions and radial functions. In literature,
ridge functions and radial functions are widely employed. Ridge functions
and radial functions are multivariate functions of the \ form $g(\mathbf{a}%
\cdot \mathbf{x})$ and $g(\left\vert \mathbf{x}-\mathbf{a}\right\vert _{e})$%
, respectively, where $\mathbf{a}\in \mathbb{R}^{d}$ is a fixed vector
(direction), $\mathbf{x}\in \mathbb{R}^{d}$ is the variable, $\mathbf{a}%
\cdot \mathbf{x}$ is the standard inner product, $\left\vert \mathbf{x-a}%
\right\vert _{e}$ is the Euclidean distance between $\mathbf{x}$ and $%
\mathbf{a}$, and $g$ is a one-variable real function.

\bigskip

\begin{center}
{\large \textbf{2. The main results}}
\end{center}

For a compact metric space $X$, and closed algebras $A_{1}\subset C(X),$ $%
A_{2}\subset C(X)$ containing the constant functions, define the following
equivalence relations in $X$:

\begin{equation*}
a\overset{R_{i}}{\sim }b\text{ if }f(a)=f(b)\text{ for all }f\in A_{i},\text{
}i=1,2.
\end{equation*}%
Then, for each $i=1,2,$ the quotient space $X_{i}=X/R_{i}$ with respect to
the relation $R_{i}$, equipped with the quotient space topology, is compact
and the natural projections $s:X\rightarrow X_{1}$ and $p:X\rightarrow X_{2}$
are continuous. Note that the quotient spaces $X_{1}$ and $X_{2}$ are not
only compact but also Hausdorff (see, e.g., \cite[p.54]{12}). In view of the
Stone-Weierstrass theorem, the algebras $A_{1}$ and $A_{2}$ have the
following set representations
\begin{eqnarray*}
A_{1} &=&\{g(s(x)):~g\in C(X_{1})\}, \\
A_{2} &=&\{h(p(x)):~h\in C(X_{2})\}.
\end{eqnarray*}

Our analysis exceedingly uses the following mathematical objects called
\textit{lightning bolts}.

\bigskip

\textbf{Definition 2.1.} (see \cite{Mar2}) \textit{A finite or infinite
ordered set $l=\{x_{1},x_{2},...\}\subset X$, where $x_{i}\neq x_{i+1}$,
with either $s(x_{1})=s(x_{2}),p(x_{2})=p(x_{3}),s(x_{3})=s(x_{4}),...$ or $%
p(x_{1})=p(x_{2}),s(x_{2})=s(x_{3}),p(x_{3})=p(x_{4}),...$ is called a
lightning bolt with respect to the algebras $A_{1}$ and $A_{2}$.}

\bigskip

For convenience, we will simply use the term \textquotedblleft
bolt\textquotedblright\ instead of the long expression \textquotedblleft
lightning bolt with respect to the algebras $A_{1}$ and $A_{2}$%
\textquotedblright . If in a finite bolt $\{x_{1},...,x_{n},x_{n+1}\}$, $%
x_{n+1}=x_{1}$ and $n$ is an even number, then the bolt $\{x_{1},...,x_{n}\}$
is said to be a \textit{closed bolt}.

Note that bolts, in the special case when $X\subset \mathbb{R}^{2}$, and $%
A_{1}=\{g(x)\}$ and $A_{2}=\{h(y)\}$, are geometrically explicit objects. In
this particular case, which is also the simplest case, a bolt is a sequence $%
\{x_{1},x_{2},...\}$ in $\mathbb{R}^{2}$ with the line segments $%
[x_{i},x_{i+1}],$ $i=1,...,n,$ perpendicular alternatively to the $x$ and $y$
axes. Such special bolts were first introduced by Diliberto and Straus in
\cite{4}. These simple objects were further utilized in many works devoted
to the approximation of multivariate functions by sums of univariate
functions (see, e.g., \cite{12} and references therein). It should be
remarked that \textquotedblleft bolt" is not the only name of these objects.
They appeared in a number of papers with several different names such as
\textit{permissible lines} (see \cite{4}), \textit{paths} (see, e.g., \cite%
{13}), \textit{trips}\ (see, e.g., \cite{Mar1}), \textit{links} (see, e.g.,
\cite{Cow}). The term \textit{bolt of lightning} is due to Arnold \cite{Arn}%
. Marshall and O'Farrell \cite{Mar2} generalized bolts to the case of two
abstract continuous function algebras and gave many central properties of
such general bolts and functionals associated with them.

In the paper \cite{AI}, \textit{extremal bolts} with respect to two algebras
were introduced.

\bigskip

\textbf{Definition 2.2.} \textit{A finite or infinite bolt \textit{$%
\{x_{1},x_{2},...\}$} is said to be extremal for a function $f\in C(X)$ if $%
f(x_{i})=(-1)^{i}\left\Vert f\right\Vert ,i=1,2,...$ or $f(x_{i})=(-1)^{i+1}%
\left\Vert f\right\Vert ,$ $i=1,2,...$}

\bigskip

In \cite{AHI}, the following Chebyshev type theorem was proved.

\bigskip

\textbf{Theorem 2.1.} \textit{Assume $X$ is a compact metric space. A
function $u_{0}\in A_{1}+A_{2}$ is a best approximation to a function $f\in
C(X)\backslash (A_{1}+A_{2})$ if and only if there exists a closed or
infinite bolt extremal for the function $f-u_{0}$.}

\bigskip

Note that this theorem was proved in \cite{AI} under quite severe
restriction that the algebras $A_{1}$ and $A_{2}$ have the $C$-property,
that is, for any $w\in C(X)$, the $\max $ and $\min $ functions
\begin{equation*}
\max_{\substack{ x\in X  \\ s(x)=a}}w(x),\text{ }\min_{\substack{ x\in X  \\ %
s(x)=a}}w(x),\text{ }\max_{\substack{ x\in X  \\ p(x)=b}}w(x),\min
_{\substack{ x\in X  \\ p(x)=b}}w(x),
\end{equation*}%
are all continuous. In the special case when $X\subset \mathbb{R}^{2}$ and $%
s,p$ are the coordinate functions, a Chebyshev type alternation theorem was
first obtained by Khavinson \cite{Kh}. In \cite{7}, similar alternation
theorems were proved for ridge functions and certain function compositions.

It should be remarked that Theorem 2.1 is an analogue of the classical
Chebyshev alternation theorem for polynomial approximation. This theorem
gives a criterion for a polynomial $P$ of degree at most $n$ to be the best
uniform approximation to a function $f\in C[a,b]$, using the oscillating
nature of the difference $f-P$. More precisely, the theorem says that $P$ is
the best approximation to $f$ if and only if there exist $n+2$ points $%
x_{i}\in \lbrack a,b]$ such that

\begin{equation*}
f(x_{i})-P(x_{i})=(-1)^{i}\max_{x\in \lbrack a,b]}\left\vert
f(x)-P(x)\right\vert ,\text{ }i=1,...,n+2.
\end{equation*}%
We refer the reader to the monograph of Natanson \cite{Nat} for a
comprehensive commentary on this theorem. For the history and variants of
the Chebyshev alternation theorem, see \cite{Bro}.

The following result is an analogue of the classical de la Vall\'{e}e
Poussin theorem.

\bigskip

\textbf{Theorem 2.2.} \textit{Assume $X$ is a compact metric space and $f\in
C(X)$. Assume $\{x_{1},x_{2},...\}$ is a closed or infinite bolt with
respect to the closed subalgebras $A_{1},A_{2}\subset C(X)$, containing
constants, and $u\in A_{1}+A_{2}$ is a function such that}

\begin{equation*}
f(x_{i})-u(x_{i})=(-1)^{i}\left\vert f(x_{i})-u(x_{i})\right\vert ,\text{
for }i=1,2,...
\end{equation*}%
\textit{Then}

\begin{equation*}
E(f)\geq \inf_{i}\left\vert f(x_{i})-u(x_{i})\right\vert .\eqno(2.1)
\end{equation*}

\bigskip

\textbf{Proof.} With each finite bolt $l=\{x_{1},...,x_{n}\}$ with respect
to $A_{1}$ and $A_{2}$, we associate the following bolt functional

\begin{equation*}
r_{l}:C(X)\rightarrow \mathbb{R},\text{ }r_{l}(F)=\frac{1}{n}%
\sum_{i=1}^{n}(-1)^{n+1}F(x_{i}).\eqno(2.2)
\end{equation*}

It is easy to verify that $r_{l}$ is a linear continuous functional with the
norm $\left\Vert r_{l}\right\Vert \leq 1$ and $\left\Vert r_{l}\right\Vert
=1 $ if and only if the set of points $x_{i}$ with odd indices $i$ does not
intersect with the set of points with even indices. Besides, if $l$ is
closed, then $r_{l}\in (A_{1}+A_{2})^{\perp }$, that is, $r_{l}(F)=0$ for
all $F\in A_{1}+A_{2}$. Here $(A_{1}+A_{2})^{\perp }$ denotes the
annihilator of the subspace $A_{1}+A_{2}\subset C(X).$ But if $l$ is not
closed, then $r_{l}$ is generally not an annihilating functional.
Nevertheless, it obeys the following important property that

\begin{equation*}
\left\vert r_{l}(F_{i})\right\vert \leq \frac{2}{n}\left\Vert
F_{i}\right\Vert ,\eqno(2.3)
\end{equation*}%
for all $F_{i}\in A_{i}$, $i=1,2$. This inequality means that for bolts $l$
with sufficiently large number of points, $r_{l}$ behaves like an
annihilating functional on each $A_{i}$, and hence on $A_{1}+A_{2}$. To see
the validity of (2.3) it is enough to recall that $F_{1}=$ $g\circ s$, $%
F_{2}=$ $h\circ p$ and consider the chain of equalities $%
g(s(x_{1}))=g(s(x_{2})),$ $g(s(x_{3}))=g(s(x_{4})),...$(or $%
g(s(x_{2}))=g(s(x_{3})),$ $g(s(x_{4}))=g(s(x_{5})),...$) for $%
F_{1}(x)=g(s(x))$ and similar equalities for $F_{2}(x)=h(p(x))$.

First assume $l=\{x_{1},x_{2},...,x_{2n}\}$ is a closed bolt. Then

\begin{equation*}
\left\vert r_{l}(f)\right\vert =\left\vert r_{l}(f-u)\right\vert =\frac{1}{2n%
}\sum_{i=1}^{2n}\left\vert f(x_{i})-u(x_{i})\right\vert \geq
\inf_{i}\left\vert f(x_{i})-u(x_{i})\right\vert .\eqno(2.4)
\end{equation*}%
Since $r_{l}\in (A_{1}+A_{2})^{\perp }$, we can write that

\begin{equation*}
\left\vert r_{l}(f)\right\vert =\left\vert r_{l}(f-F)\right\vert \leq
\left\Vert f-F\right\Vert \eqno(2.5)
\end{equation*}%
for any $F\in A_{1}+A_{2}$. It follows from (2.5) and arbitrariness of $F$
that

\begin{equation*}
\left\vert r_{l}(f)\right\vert \leq E(f).\eqno(2.6)
\end{equation*}%
Inequalities (2.4) and (2.6) together yield (2.1).

Assume now $l=\{x_{1},x_{2},...\}$ is an infinite bolt. In this case we will
use the famous Banach-Alaoglu theorem on weak$^{\text{*}}$ sequential
compactness of the closed unit ball in $E^{\ast }$ for a separable Banach
space $E$ (see, e.g., Rudin \cite[p. 66]{Rud}). Note that the space $C(X)$
is separable, since $X$ is a compact metric space. Hence the closed unit
ball $B$ in $C^{\ast }(X)$ is sequentially compact. It follows that any
sequence in $B$ has a weak$^{\text{*}}$ convergent subsequence.

Now we proceed as follows. From $l$ we form the finite bolts $%
l_{k}=\{x_{1},...,x_{k}\},$ $k=1,2,...$, and consider the bolt functionals $%
r_{l_{k}}$. For the ease of notation, let us put $r_{k}=r_{l_{k}}.$ The
sequence $\{r_{_{k}}\}_{k=1}^{\infty }$ is contained in the closed unit ball
$B$ of the dual space $C^{\ast }(X).$ By the Banach-Alaoglu theorem, this
sequence has a weak$^{\text{*}}$ convergent subsequence. Without loss of
generality we may assume that the sequence $\{r_{_{k}}\}_{k=1}^{\infty }$
itself weak$^{\text{*}}$ converges to some element $r^{\ast }\in B$. From
(2.3) it follows that $r^{\ast }(F_{1}+F_{2})=0,$ for any $F_{i}\in A_{i},$ $%
i=1,2$. That is, $r^{\ast }\in B\cap (A_{1}+A_{2})^{\perp }$.

For each functional $r_{k}$ we can write that

\begin{equation*}
\left\vert r_{k}(f-u)\right\vert =\frac{1}{k}\sum_{i=1}^{k}\left\vert
f(x_{i})-u(x_{i})\right\vert \geq \inf_{i}\left\vert
f(x_{i})-u(x_{i})\right\vert ,\text{ }k=1,2,...\eqno(2.7)
\end{equation*}%
From (2.7) it follows that

\begin{equation*}
\left\vert r^{\ast }(f-u)\right\vert \geq \inf_{i}\left\vert
f(x_{i})-u(x_{i})\right\vert .\eqno(2.8)
\end{equation*}%
Since $r^{\ast }$ annihilates all elements of $A_{1}+A_{2}$ and $\left\Vert
r^{\ast }\right\Vert \leq 1$, it follows that

\begin{equation*}
\left\vert r^{\ast }(f-F)\right\vert \leq E(f),\eqno(2.9)
\end{equation*}%
for all $F\in A_{1}+A_{2}$. Taking now $F=u$ in (2.9) and considering (2.8),
we obtain the validity of (2.1). The theorem has been proved.

\bigskip

In the sequel, we obtain a formula for the approximation error $E(f)$. The
well-known duality relation says that

\begin{equation*}
E(f)=\underset{\left\Vert G\right\Vert \leq 1}{\underset{G\in
(A_{1}+A_{2})^{\perp }}{\sup }}\left\vert G(f)\right\vert ,\eqno(2.10)
\end{equation*}%
The $\sup $ in (2.10) is attained by some functional $G$ with the norm $%
\left\Vert G\right\Vert =1.$ We are interested in the problem: is it
possible to replace in (2.10) the annihilator class $(A_{1}+A_{2})^{\perp }$
by some subclass of it consisting of functionals of simple structure? In the
following we see that under a reasonable assumption, the functionals $r_{l}$
generated by closed bolts $l$ suffice. Namely, we assume proximinality of
the subspace $A_{1}+A_{2}$ in $C(X)$. Let $E$ be a normed linear space and $%
F $ be its subspace. We say that $F$ is proximinal in $E$ if for any element
$e\in E$ there exists at least one element $f_{0}\in F$ such that

\begin{equation*}
\left\Vert e-f_{0}\right\Vert =\inf_{f\in F}\left\Vert e-f\right\Vert .
\end{equation*}%
Proximinality of the sum $A_{1}+A_{2}$ in some special cases of algebras $%
A_{i}$ was investigated in \cite{Gar1,Gar2,6}. We refer the reader to Alimov
and Tsar'kov's book \cite{A} for the general theory of proximinal sets and
for the study of their various properties from the viewpoint of geometric
approximation theory.

The following theorem is valid.

\bigskip

\textbf{Theorem 2.3.} \textit{Assume $X$ is a compact metric space, $%
A_{1},A_{2}\subset C(X)$ are closed algebras containing constants and $%
A_{1}+A_{2}$ is proximinal in $C(X)$. Then the approximation error $E(f)$ of
a function $f\in C(X)\backslash (A_{1}+A_{2})$ obeys the the following
duality formula}

\begin{equation*}
E(f)=\sup_{l\subset X}\left\vert r_{l}(f)\right\vert ,\eqno(2.11)
\end{equation*}%
\textit{where the $\sup $ is taken over all closed bolts.}

\bigskip

\textbf{Proof.} First note that on the basis of (2.6) the one-side
inequality is always true:

\begin{equation*}
\sup_{l\subset X}\left\vert r_{l}(f)\right\vert \leq E(f).\eqno(2.12)
\end{equation*}

By assumption, $f$ has a best approximation in $A_{1}+A_{2}$. Denote this
function by $u_{0}$. Let us concentrate on extremal bolts for the function $%
f_{1}=f-u_{0}$. Theorem 2.1 says that regarding such bolts there may be only
two cases.

\bigskip

\textbf{Case 1.} There exists a closed bolt $l_{0}=(x_{1},...,x_{2n})$
extremal for the function $f_{1}$.

In this case, based on closedness of $l_{0}$ and Definition 2.2, we can
write that

\begin{equation*}
\left\vert r_{l_{0}}(f)\right\vert =\left\vert r_{l_{0}}(f-u_{0})\right\vert
=\left\Vert f-u_{0}\right\Vert =E(f).\eqno(2.13)
\end{equation*}%
(2.12) and (2.13) together yield (2.11).

\bigskip

\textbf{Case 2.} There exists an infinite bolt $l=\{x_{1},x_{2},...\}$
extremal for $f_{1}$.

Since $A_{1}+A_{2}$ is proximinal in $C(X)$, clearly $A_{1}+A_{2}$ is closed
in $C(X)$. On the other hand, $A_{1}+A_{2}$ is closed in $C(X)$ if and only
if the lengths (number of points) of irreducible bolts are uniformly bounded
(see \cite{Med}). Note that a bolt $\{x_{1},x_{2},...,x_{n}\}$ is called
\textit{irreducible} if there is not a bolt connecting $x_{1}$ and $x_{n}$
with the length less than $n$.

Let a positive integer $N$ bound the lengths of all irreducible bolts in $X$%
. Then the finite extremal bolts $\{x_{1},x_{2},...,x_{n}\}\subset l$, $%
n=N+1,N+2,...$, or subbolts of them, must be made closed by adding not more
than $N$ points. Without loss of generality we may assume that these bolts
themselves can be made closed. That is, for each finite extremal bolt $%
l_{n}=\{x_{1},x_{2},...,x_{n}\}$, $n>N$, there exists a closed bolt $%
L_{n}=(x_{1},x_{2},...,x_{n},y_{n+1},...,y_{n+m_{n}})$, where $m_{n}\leq N$.
The functional $r_{L_{n}}$ obeys the inequalities

\begin{equation*}
\left\vert r_{L_{n}}(f)\right\vert =\left\vert r_{L_{n}}(f-u_{0})\right\vert
\leq \frac{n\left\Vert f-u_{0}\right\Vert +m_{n}\left\Vert
f-u_{0}\right\Vert }{n+m_{n}}=\left\Vert f-u_{0}\right\Vert \eqno(2.14)
\end{equation*}%
and
\begin{equation*}
\left\vert r_{L_{n}}(f)\right\vert \geq \frac{n\left\Vert f-u_{0}\right\Vert
-m_{n}\left\Vert f-u_{0}\right\Vert }{n+m_{n}}=\frac{n-m_{n}}{n+m_{n}}%
\left\Vert f-u_{0}\right\Vert .\eqno(2.15)
\end{equation*}%
We obtain from (2.14) and (2.15) that
\begin{equation*}
\sup_{L_{n}}\left\vert r_{L_{n}}(f)\right\vert =\left\Vert
f-u_{0}\right\Vert =E(f).\eqno(2.16)
\end{equation*}%
The inequality (2.16) together with (2.12) yield that

\begin{equation*}
E(f)=\sup_{l\subset X}\left\vert r_{l}(f)\right\vert ,
\end{equation*}%
where the $\sup $ is taken over all closed bolts. The theorem has been
proved.

\bigskip

Theorem 2.3 for the case in which $X$ is a rectangle in $\mathbb{R}^{2}$
with sides parallel to the coordinate axis and the quotient mappings $s,p$
are the coordinate projections was proved by Diliberto and Straus \cite{4},
and independently by Ofman \cite{Ofm}. Note that by Kolmogorov's theorem for
such sets $X$ the set $A_{1}+A_{2}=\left\{ g(x)+h(y)\right\} $ is always
proximinal in $C(X)$ (see \cite{Ofm}). The same result for arbitrary compact
sets $X\subset \mathbb{R}^{2}$ but with the assumption that $\left\{
g(x)+h(y)\right\} $ is proximinal was proved by Khavinson \cite{Kh}. In this
special case, Khavinson also proved Theorem 2.2.

A duality relation for the approximation error in the abstract algebraic
case was obtained by Marshall and O'Farrell \cite{Mar2}. In \cite{Mar2},
together with other results, they proved the validity of the formula

\begin{equation*}
E(f)=\max_{l\subset X}\left\{ \underset{n\rightarrow \infty }{\lim \sup }%
\left\vert r_{l_{n}}(f)\right\vert \right\} ,\eqno(2.17)
\end{equation*}%
where the $\max $ is taken over all infinite bolts $l=\{x_{1},...,x_{n},...%
\} $ and $l_{n}=\{x_{1},...,x_{n}\}$. Note that formula (2.17) holds for a
compact Hausdorff space $X$ and without the proximinality hypothesis, but
with the proviso that the algebras $A_{1}$ and $A_{2}$ together separate
points of $X$, that is, for any $x_{1}\neq x_{2}$ in $X$ there are functions
$G\in A_{1}$ and $H\in A_{2}$ such that $(G(x_{1}),H(x_{1}))\neq
(G(x_{2}),H(x_{2}))$. For example, this condition holds in the above case $%
A_{1}+A_{2}=\left\{ g(x)+h(y)\right\} $. But it does not hold in the case of
ridge function algebras $A_{1}=\{g(\mathbf{a}\cdot \mathbf{x})\}$ and $%
A_{2}=\{h(\mathbf{b}\cdot \mathbf{x})\}$, where $\mathbf{a}$ and $\mathbf{b}$
are fixed vectors in $\mathbb{R}^{d}$, but $g$ and $h$ are variable
functions. Also note that the set of infinite bolts in (2.17) contains
closed bolts as well. This is because a closed bolt $l=\{x_{1},...,x_{2n}\}$
can be written as an infinite bolt $l=%
\{x_{1},...,x_{2n},x_{1},...,x_{2n},...\}$ by infinite repetition of its
points. In this case,

\begin{equation*}
\underset{n\rightarrow \infty }{\lim \sup }\left\vert
r_{l_{n}}(f)\right\vert =\left\vert r_{l}(f)\right\vert ,
\end{equation*}%
and hence%
\begin{equation*}
\sup_{l_{\text{closed}}\subset X}\left\vert r_{l}(f)\right\vert \leq
\max_{l_{\text{infinite}}\subset X}\left\{ \underset{n\rightarrow \infty }{%
\lim \sup }\left\vert r_{l_{n}}(f)\right\vert \right\} \leq E(f).\eqno(2.18)
\end{equation*}%
The $\sup $ and $\max $ in (2.18) are taken over all closed bolts and all
infinite bolts, respectively. The right inequality in (2.18) follows from
the fact that for any infinite bolt $l=\{x_{1},...,x_{n},...\}$ and finite
bolts $l_{n}=\{x_{1},...,x_{n}\}$, the functionals $r_{l_{n}}$ are in the
closed unit ball of $C^{\ast }(X)$ and a limit of any weak$^{\text{*}}$
convergent subsequence of $\left\{ r_{l_{n}}\right\} _{n=1}^{\infty }$
belongs to $(A_{1}+A_{2})^{\perp }$. Inequality (2.18) shows that in case if
$A_{1}+A_{2}$ is proximinal in $C(X)$, the formula (2.11) is more efficient
than (2.18), since in (2.11)\ it is not needed to consider infinite bolts
different from closed bolts. The $\sup $ in (2.11) is attained if there is a
closed extremal bolt for the function $f$ (see Case 1 above).

\end{document}